# About one method for constructing Hermite trigonometric polynomials

Denysiuk V.P., Dr. of Phys-Math. Sciences, Professor, Kyiv, Ukraine

National Aviation University

kvomden@nau.edu.ua

## Annotation

The method of constructing Hermite trigonometric polynomials, which interpolate the values of a certain periodic function and its derivatives up to (including ) the $p$-th ( $p = 1, 2, \ldots$ ) order in nodes of a uniform grid, is considered. The proposed method is based on the properties of the periodicity of trigonometric functions and is reduced to the solving of only two systems of the linear algebraic equations of the $p+1$--th order; solutions to these systems may be obtained in advance. When implementing this method, it is necessary to calculate the coefficients of interpolation trigonometric polynomials that interpolate the values of the function itself and the values of its derivatives in the nodes of a uniform grid; well-known algorithms of fast Fourier transform can be applied for this purpose. The examples of construction of Hermite trigonometric polynomials for $p = 1, 2 \ldots$ are given. The proposed method may be recommended for wide practical application.

## Introduction

In many scientific and technological problems the statement of problem of interpolation of functions in such formulation often arises. Let on a segment $[a,b]$ a certain grid $a = x_0 < x_1 < \ldots < x_N = b$ be given and let the values of the function itself $f(x_i) = f_i$ ( $i = 0, 1, \ldots, N$ ), and the values of the derivatives, $f^{(k)}(x_i) = f_i^{(k)}$, ( $k = 1, 2, \ldots, p$ ) of this function up to (including ) $p$-th ( $p = 1, 2, \ldots$ ) order in the grid nodes be known. It is necessary to construct a function $H(x)$ of a certain class, which depends on the parameters $(N+1)(p+1)$ and satisfies the conditions

$$H^{(k)}(x_i) = f_i^{(k)}, \qquad (i = 0, 1, \ldots, N; k = 0, 1, \ldots, p). \tag{1}$$

In the case where the function $H(x) \in P$, where P is the class of algebraic polynomials, then the polynomial satisfying these conditions is called the Hermite interpolation polynomial.

In many cases, a uniform grid is In many cases, a uniform grid is given on a segment $[0, 2\pi)$ and a function $H(x)$ is searched in a class of trigonometric polynomials T; a polynomial $T(x)$, which is the solution of such a problem, is called an Hermite trigonometric polynomial. The problem of constructing an

Hermite trigonometric polynomial satisfying the conditions (1) will be considered further.

## Analysis of researches and publications

There are different approaches to the construction of Hermite trigonometric polynomials. One of them, for example, is based on the modification of methods of polynomial Hermite polynomial interpolation. Note that the problem of constructing a power Hermite polynomial in general terms is considered in detail, for example, in [1]. Partial cases of constructing a power polynomial are considered in [2]. The construction of the same Hermite trigonometric polynomials and the problems associated with this approach are considered, for example, in [3]; an adequate bibliography on this subject is also provided. This approach and its development are also considered in papers [4], [5].

Another approach to constructing Hermite trigonometric polynomials is based on the application of fundamental trigonometric functions. Such approach is considered, for example, in [6]; it is worth pointing out that in this work, in our opinion, fundamental trigonometric functions are constructed most naturally. Other variants of this approach, which differ in methods of constructing fundamental functions, are discussed in papers [7], [8], [9]. The disadvantage of this approach is the complexity of constructing fundamental trigonometric functions of higher orders; the consequence is that the constructed Hermite polynomials interpolate the function itself and its derivative only in the first order.

## The purpose of the work.

The purpose of the work is to develop a constructive method for constructing Hermite trigonometric polynomial, which interpolates the function and its derivatives of order $p$, ( $p = 1, 2, \ldots$ ) in nodes of uniform grids on the interval $[0, 2\pi)$; this method is based on the property of the periodicity of trigonometric functions.

## The main part.

When interpolating functions, the case is often encountered when the nodes of a uniform grid $\Delta_N^t$ given on a segment $[0, 2\pi]$, we know not only the values of the function $f(t)$, but also the values of the derivatives of this function. In the case where interpolation is carried out by algebraic polynomials, the Hermite polynomials are used, which interpolate in the nodes of the grid $\Delta_N^t$ both the values of the function itself and the values of its derivatives. If the interpolation of a function is supposed to be carried out by trigonometric polynomials, then the question arises about the construction of trigonometric polynomials, which are analogues of Hermite algebraic polynomials; in the following, these trigonometric polynomials will be called Hermite trigonometric polynomial. Let's consider the construction of Hermite trigonometric polynomials in more detail.

Let the grid be given to $[0,2\pi)$ $\Delta_N^{(I)} = \left\{ t_j^{(I)} \right\}_{j=1}^{N}$, $(I = 0,1; \ N = 2n+1, n = 1,2,...;$

$t_j^{(0)} = \frac{2\pi}{N}(j-1)$; $t_j^{(1)} = \frac{\pi}{N}(2j-1)$), and let in the nodes of this grid the values of the

function $f(t)$ and its derivatives up to p-th the order inclusive $(p = 1,2,...)$ be given.

We denote the values of the function and its derivatives in the nodes of the grid $\Delta_N^{(I)}$

by $f_j^{(I,0)}, f_j^{(I,1)},..., f_j^{(I,p)}$, $(j = 1,2,...,N)$. With respect to these derivatives, we will

assume that these derivatives are centered on the grid $\Delta_N^{(I)}$, that is, the condition

$$\sum_{j=1}^{N} f_j^{(I,i)} = 0, \ \ (i = 1,2,...,p)$$

is fulfilled.

In the future, we will have to distinguish cases of even and odd values of the parameter p.

The Hermite trigonometric polynomial that interpolates the function and its derivatives to the p-th order of the order inclusively ( p is odd, $p = 2q-1, q = 1,2,...$) in the nodes of the grid $\Delta_N^I, (N = 2n+1, n = 1,2,...)$, will be sought in the form

$$T_{(p+1)n}^{(I)}(t) = \frac{a_0^{(I)}}{2} +$$

(1)

$$+ \sum_{k=1}^{n} \left\{ a_k^{(I)} \cos kt + \sum_{i=1}^{q-1} \left[ a_{iN-k}^{(I)} \cos(iN-k)t + a_{iN+k}^{(I)} \cos(iN+k)t \right] + a_{qN-k}^{(I)} \cos(qN-k)t \right\} +$$

$$+ \sum_{k=1}^{n} \left\{ b_k^{(I)} \sin kt + \sum_{i=1}^{q-1} \left[ b_{iN-k}^{(I)} \sin(iN-k)t + b_{iN+k}^{(I)} \sin(iN+k)t \right] + b_{qN-k}^{(I)} \sin(qN-k)t \right\}.$$

To find the coefficients $a_j^{(I)}, b_j^{(I)}$, we need to do next steps. We compute the

values of a polynomial $T_{(p+1)n}^{(I)}(t)$ in the nodes of the grid $\Delta_N^{(I)}$ and equate them to

the values of the trigonometric polynomial $T_{0,n}^{(I)}(t)$ that interpolates the values of

the function in the nodes of the same grid. Then, by differentiating the polynomial

$T_{(p+1)n}^{(I)}(t)$ $k$ times $(k = 1,2,...,p)$ and equating the values of the derivative of k-order

in the nodes of the grid $\Delta_N^{(I)}$, to the values of the interpolation trigonometric polynomial $T_{k,n}^{(I)}(t)$, interpolating the values of the derivative of the same order of function $f(t)$ in the nodes of the grid $\Delta_N^{(I)}$, we obtain 2n systems of linear algebraic equations of $p+1$ -th order of the desired coefficients of the polynomial $T_{(p+1)n}^{(I)}(t)$.

The proposed approach to constructing Hermite trigonometric polynomial is based on the fact that the values of trigonometric functions $\cos(iN \pm k)t$ and

$\sin(iN \pm k)t$, in the nodes of the grid $\Delta_N^{(I)}$, depend only on the numbers of the nodes

of this grid j and the indices k; that is, equalities $\cos(iN \pm k)t_j^{(I)} = \left| \cos kt_j^{(I)} \right|$,

$\sin(iN \pm k)t_j^{(I)} = \left| \sin kt_j^{(I)} \right|$ are fulfilled. The use of this fact makes it possible to

simplify the problem of finding the coefficients $a_j, b_j$, of a polynomial $T_{(q+1)n}^{(I)}(t)$, by reducing it to the problem of calculating the coefficients of interpolation trigonometric polynomials and then solving the N-1 systems of algebraic equations of dimension p+1.

We illustrate the proposed approach on the example of constructing the Hermite trigonometric polynomial that interpolates the function $f(t)$ and its derivative of the first order in the nodes of the grid $\Delta_N^{(0)}$. In this case $p=1$, $q=1$; the internal sum disappears, and the Hermite trigonometric polynomial will be sought in the form

$$T_{2n}^{(0)}(t) = \frac{a_0^{(0)}}{2} + \sum_{k=1}^{n} \left[ a_k^{(0)} \cos kt + a_{N-k}^{(0)} \cos(N-k)t + \right. \tag{2}$$
$$\left. + b_k^{(0)} \sin kt + b_{N-k}^{(0)} \sin(N-k)t \right].$$

Calculate the values of the polynomial $T_{2n}^{(0)}(t)$ in the nodes of the grid $\Delta_N^0$. Taking account of that $t_i^{(0)} = \frac{2\pi}{N}(j-1)$ and equality

$$\cos(N-k)t_j^{(0)} = \cos(N-k)\frac{2\pi}{N}(j-1) = \cos kt_j^{(0)},$$
$$\sin(N-k)t_j^{(0)} = \sin(N-k)\frac{2\pi}{N}(j-1) = -\sin kt_j^{(0)},$$

we have

$$T_{2n}^{(0)}(t_j^{(0)}) = \frac{a_0^{(0)}}{2} + \sum_{k=1}^{n} \left[ (a_k^{(0)} + a_{N-k}^{(0)}) \cos kt_j^{(0)} + (b_k^{(0)} - b_{N-k}^{(0)}) \sin kt_j^{(0)} \right].$$

Hence it follows that the coefficients $a_k^{(0)}, a_{N-k}^{(0)}, b_k^{(0)}, b_{N-k}^{(0)}$, satisfy the systems of equations

$$a_0^{(0)} = A_{0,0}^{(0)};$$
$$\begin{cases} a_k^{(0)} - a_{N-k}^{(0)} = A_{0,k}^{(0)} \\ b_k - b_{N-k} = B_{0,k}^{(0)}: \qquad k=1,2,...,n. \end{cases} \tag{3}$$

where $A_{0,0}^{(0)}$, $A_{0,k}^{(0)}$ and $B_{0,k}^{(0)}$ are coefficients of the trigonometric polynomial that interpolates the function f(t) in the nodes of the grid $\Delta_N^{(0)}$.

Differentiating the polynomial $T_{2n}^{(0)}(t)$, we have

$$T_{2n}^{(0)'}(t) = \sum_{k=1}^{n} \left[ (-k)a_k^{(0)} \sin kt - (N-k)a_{N-k}^{(0)} \sin(N-k)t + \right.$$
$$\left. + kb_k^{(0)} \cos kt + (N-k)b_{N-k}^{(0)} \cos(N-k)t \right].$$

Calculating values of the polynomial $T_{2n}^{(0)'}(t)$ in the nodes of the grid $\Delta_N^0$, we obtain

$$T_{2n}^{(0)'}(t_j^{(0)}) = \sum_{k=1}^{n} \left\{ \left[ kb_k^{(0)} + (N-k)b_{N-k}^{(0)} \right] \cos kt_j^{(0)} + \right.$$
$$\left. + \left[ (-k)a_k^{(0)} + (N-k)a_{N-k}^{(0)} \right] \sin kt_j^{(0)} \right\}.$$

Denoting by $A_{1,k}^{(0)}, B_{1,k}^{(0)}$ coefficients of the trigonometric polynomial interpolating the values of the derivative $f'(t)$ in the nodes of the grid $\Delta_N^{(0)}$, and taking into account that $A_{1,0}^{(0)} = 0$ we, as before, obtain systems of equations

$$\begin{cases} kb_k^{(0)} + (N-k)b_{N-k}^{(0)} = A_{1,k}^{(0)} \\ (-k)a_k^{(0)} + (N-k)a_{N-k}^{(0)} = B_{1,k}^{(0)}; \qquad k=1,2,...,n. \end{cases} \qquad (4)$$

Systems of equations (3), (4) may be presented in the form

$$\begin{cases} a_k^{(0)} + a_{N-k}^{(0)} = A_{0,k}^{(0)} \\ (-k)a_k^{(0)} + (N-k)a_{N-k}^{(0)} = B_{1,k}^{(0)}; \end{cases} \qquad \begin{cases} b_k^{(0)} - b_{N-k}^{(0)} = B_{0,k}^{(0)} \\ kb_k^{(0)} + (N-k)b_{N-k}^{(0)} = A_{1,k}^{(0)}; \end{cases}$$

$$k = 1,2,...,n; \qquad\qquad\qquad k = 1,2,...,n.$$

Solving these systems of equations, we have

$$a_k^{(0)} = \frac{1}{N}\left[(N-k)A_{0,k}^{(0)} - B_{1,k}^{(0)}\right]; \qquad a_{N-k}^{(0)} = \frac{1}{N}\left[kA_{0,k}^{(0)} + B_{1,k}^{(0)}\right];$$

$$b_k^{(0)} = \frac{1}{N}\left[(N-k)B_{0,k}^{(0)} + A_{1,k}^{(0)}\right]; \qquad b_{N-k}^{(0)} = \frac{1}{N}\left[-kB_{0,k}^{(0)} + A_{1,k}^{(0)}\right];$$

$$k = 1,2,...,n.$$

So, the Hermite trigonometric polynomial that interpolates the function f(t) and its derivative in the nodes of the grid $\Delta_N^{(0)}$ has the form

$$T_{2n}^{(0)}(t) = \frac{A_{0,0}^{(0)}}{2} + \frac{1}{N}\sum_{k=1}^{n}\left\{\left[(N-k)A_{0,k}^{(0)} - B_{1,k}^{(0)}\right]\cos kt + \left[kA_{0,k}^{(0)} + B_{1,k}^{(0)}\right]\cos(N-k)t + \right.$$
$$\left. + \left[(N-k)B_{0,k}^{(0)} + A_{1,k}^{(0)}\right]\sin kt + \left[-kB_{0,k}^{(0)} + A_{1,k}^{(0)}\right]\sin(N-k)t\right\}. \qquad (5)$$

Similarly, it is easy to obtain a Hermite trigonometric polynomial that interpolates the function f(t) and its derivative in the nodes of a grid $\Delta_N^1$, that has the form

$$T_{2n}^{(1)}(t) = \frac{A_{0,0}^{(1)}}{2} + \frac{1}{N}\sum_{k=1}^{n}\left\{\left[(N-k)A_{0,k}^{(1)} - B_{1,k}^{(1)}\right]\cos kt - \left[kA_{0,k}^{(1)} + B_{1,k}^{(1)}\right]\cos(N-k)t + \right.$$
$$\left. + \left[(N-k)B_{0,k}^{(1)} + A_{1,k}^{(1)}\right]\sin kt + \left[kB_{0,k}^{(1)} - A_{1,k}^{(1)}\right]\sin(N-k)t\right\}. \qquad (6)$$

where the coefficients $A_{0,0}^{(1)}, A_{0,k}^{(1)} \; B_{0,k}^{(1)}, \; A_{1,k}^{(1)}, \; B_{1,k}^{(1)}$, are the coefficients of interpolation trigonometric polynomials that interpolates the function and its derivative respectively in the nodes of the grid $\Delta_N^{(1)}$.

In the case when p is even, $p = 2q, \quad q = 1,2,...,$ Hermite trigonometric polynomial that interpolates the function f(t) and its derivatives up to p-th order inclusive in the nodes of the grid $\Delta_N^{(l)}$, will be sought in the form

$$T_{(q+1)n}^{(l)}(t) = \frac{a_0^{(l)}}{2} +$$
$$+ \sum_{k=1}^{n}\left\{a_k^{(l)}\cos kt + \sum_{i=1}^{q}\left[a_{iN-k}^{(l)}\cos(iN-k)t + a_{iN+k}^{(l)}\cos(iN+k)t\right]\right\} +$$
$$+ \sum_{k=1}^{n}\left\{b_k^{(l)}\sin kt + \sum_{i=1}^{q}\left[b_{iN-k}^{(l)}\sin(iN-k)t + b_{iN+k}^{(l)}\sin(iN+k)t\right]\right\}. \qquad (7)$$

Calculating the values of a polynomial $T^{(I)}_{(p+1)n}(t)$ and its derivatives in the nodes of the grid $\Delta^{(I)}_N$ and equating them with the corresponding values of trigonometric polynomials that interpolate the function and its derivatives in the nodes of this grid, it is easy to obtain systems of equations for the desired coefficients $a^{(I)}_k, b^{(I)}_k$ of the polynomial $T^{(I)}_{(p+1)n}(t)$. We restrict ourselves to the fact that we give the trigonometric Hermite polynomials that interpolate the function f(t) and its derivatives up to the second order inclusive, in the nodes of the grids $\Delta^{(0)}_N$ and $\Delta^{(1)}_N$. Such polynomials have the form

$$T^{(I)}_{3n}(t) = \frac{A^{(I)}_{0,0}}{2} +$$
$$+ \sum_{k=1}^{n} \Big[ a^{(I)}_k \cos kt + a^{(I)}_{N-k} \cos(N-k)t + a^{(I)}_{N+k} \cos(N+k)t +$$
$$+ b^{(I)}_k \sin kt + b^{(I)}_{N-k} \sin(N-k)t + b^{(I)}_{N+k} \sin(N+k)t \Big], \qquad (8)$$

where for $I=0$ we have

$$a^{(0)}_k = \frac{1}{N^2} \Big[ A^{(0)}_{0,k}(N^2 - k^2) - 2kB^{(0)}_{1,k} + A^{(0)}_{2,k} \Big];$$

$$b^{(0)}_k = \frac{1}{N^2} \Big[ B^{(0)}_{0,k}(N^2 - k^2) + 2kA^{(0)}_{1,k} + B^{(0)}_{2,k} \Big];$$

$$a^{(0)}_{N-k} = \frac{1}{2N^2} \Big[ A^{(0)}_{0,k}(kN + k^2) + B^{(0)}_{1,k}(N + 2k) - A^{(0)}_{2,k} \Big];$$

$$b^{(0)}_{N-k} = \frac{1}{2N^2} \Big[ -B^{(0)}_{0,k}(kN + k^2) + A^{(0)}_{1,k}(N + 2k) + B^{(0)}_{2,k} \Big];$$

$$a^{(0)}_{N+k} = \frac{1}{2N^2} \Big[ A^{(0)}_{0,k}(k^2 - kN) + B^{(0)}_{1,k}(2k - N) - A^{(0)}_{2,k} \Big];$$

$$b^{(0)}_{N+k} = \frac{1}{2N^2} \Big[ B^{(0)}_{0,k}(k^2 - kN) + A^{(0)}_{1,k}(N - 2k) - B^{(0)}_{2,k} \Big];$$

$$k = 1, 2, ..., n.$$

Accordingly, for $I=1$ we have

$$a^{(1)}_k = \frac{1}{N^2} \Big[ A^{(1)}_{0,k}(N^2 - k^2) - 2kB^{(1)}_{1,k} + A^{(1)}_{2,k} \Big];$$

$$b^{(1)}_k = \frac{1}{N^2} \Big[ B^{(1)}_{0,k}(N^2 - k^2) + 2kA^{(1)}_{1,k} + B^{(1)}_{2,k} \Big];$$

$$a^{(1)}_{N-k} = \frac{1}{2N^2} \Big[ -A^{(1)}_{0,k}(kN + k^2) - B^{(1)}_{1,k}(N + 2k) + A^{(1)}_{2,k} \Big];$$

$$b^{(1)}_{N-k} = \frac{1}{2N^2} \Big[ B^{(1)}_{0,k}(kN + k^2) - A^{(1)}_{1,k}(N + 2k) - B^{(1)}_{2,k} \Big];$$

$$a^{(1)}_{N+k} = \frac{1}{2N^2} \Big[ A^{(1)}_{0,k}(kN - k^2) + B^{(1)}_{1,k}(N - 2k) + A^{(1)}_{2,k} \Big];$$

$$b^{(1)}_{N+k} = \frac{1}{2N^2} \Big[ B^{(1)}_{0,k}(kN - k^2) - A^{(1)}_{1,k}(N - 2k) + B^{(1)}_{2,k} \Big];$$

$$k = 1, 2, ..., n;$$

where $A^{(1)}_{0,0}$, $A^{(1)}_{0,k}$, $B^{(1)}_{0,k}$ are the coefficients of the trigonometric polynomial that interpolates the function f(t) in the nodes of the grid $\Delta^{(1)}_N$, and $A^{(1)}_{1,k}, B^{(1)}_{1,k}$ and $A^{(1)}_{2,k}, B^{(1)}_{2,k}$ are the coefficients of trigonometric polynomials that respectively interpolate the derivatives of the 1st and 2nd orders of functions f(t) in the nodes of this grid.

Finally, it is advisable to make the following remark. When constructing trigonometric Hermite polynomials, we assumed that the values of the derivatives on the grid $\Delta_N^{(I)}$ are centered. However, this restriction is easy to remove, for example, in such way. Consider the differential analogue of the Kronecker symbol

$$\frac{d^j}{dx^j}H_i = \begin{cases} 1, \ i = j; \\ 0, i \neq j, \end{cases} \quad i, j = 0, 1, ..., p.$$

Taking into account this, the Hermite trigonometric polynomial, for example, for odd values p, will be sought in the form

$$T_{(q+1)n}^{(I)}(t) = \frac{a_{0,0}^{(I)}}{2}H_0 + \frac{a_{1,0}^{(I)}}{2}H_1 + ... + \frac{a_{p,0}^{(I)}}{2}H_p +$$

$$+ \sum_{k=1}^{n}\left\{ a_k^{(I)}\cos kt + \sum_{i=1}^{q-1}\left[ a_{iN-k}^{(I)}\cos(iN-k)t + a_{iN+k}^{(I)}\cos(iN+k)t \right] + a_{qN-k}^{(I)}\cos(qN-k)t \right\} +$$

$$+ \sum_{k=1}^{n}\left\{ b_k^{(I)}\sin kt + \sum_{i=1}^{q-1}\left[ b_{iN-k}^{(I)}\sin(iN-k)t + b_{iN+k}^{(I)}\sin(iN+k)t \right] + b_{qN-k}^{(I)}\sin(qN-k)t \right\},$$

where

$$a_{k,0}^{(I)} = \frac{1}{N}\sum_{j=1}^{N}f_j^{(I,k)} \ , \ (k = 0, 1, ..., p) \ .$$

In the future we will assume that when calculating the values of the Hermite trigonometric polynomial we apply an operator $d^0\big/dx^0$ to it; herewith

$$\frac{d^0}{dx^0}H_0 = 1, \quad \frac{d^0}{dx^0}H_1 = 0, \quad \frac{d^0}{dx^0}H_2 = 0, \ \ldots, \ \frac{d^0}{dx^0}H_p = 0$$

and we have formula (2).

When calculating the same values of the derivative i- order ($i = 1, 2, ..., p$) we have

$$\frac{d^i}{dx^i}H_0 = 1, \quad \frac{d^i}{dx^i}H_1 = 0, \ \ldots, \ \frac{d^i}{dx^i}H_i = 1, \ \ldots, \ \frac{d^0}{dx^0}H_p = 0$$

and we obtain a formula for the derivative with a free term.

Thus, the proposed approach is reduced to the solving of two systems of algebraic equations of dimension p+1; in addition, well-known algorithms for fast Fourier transform may be applied for calculation of the right-hand sides of the equations of these systems.

## CONCLUSIONS

1. A method is proposed for constructing Hermite trigonometric polynomials that interpolate the function itself and its centered derivatives up to the p-th p=1,2,… order inclusive, in the nodes of the grid $\Delta_N^{(I)}$, ($N = 2n+1$, $n = 1, 2, ...$), for implementation of which it is necessary to solve two systems of linear algebraic equations, each of which is a system p+1-th order.

2. The solving of systems of linear algebraic equations necessary for the implementation of the proposed method can be carried out in advance; for practical application in each case it is sufficient to substitute the coefficients of interpolation trigonometric polynomials to the solutions found.

3. In calculating the coefficients of interpolation trigonometric polynomials necessary for the construction of Hermite trigonometric polynomial using the proposed method it is advisable to use well-known algorithms for the fast Fourier transform .

4. Taking into account the above mentioned, the proposed method for constructing Hermite trigonometric polynomials, in our opinion, can be recommended for wide use in practice.

## List of References